\documentclass{article}
\usepackage[english]{babel}
\usepackage{amsthm}
\usepackage{amsfonts}
\usepackage{amssymb}
\usepackage{amsmath}
\usepackage{graphicx}
\usepackage[dvipsnames]{xcolor}
\usepackage{float}
\usepackage{xcolor}
\usepackage{indentfirst}
 \usepackage{theoremref}
\usepackage{float}

\newtheorem*{theorem*}{Theorem}
\newtheorem{theorem}{Theorem}
\newtheorem{corollary}{Corollary}[theorem]
\newtheorem{lemma}[theorem]{Lemma}

\newtheoremstyle{named}{}{}{\itshape}{}{\bfseries}{.}{.5em}{\thmnote{#3's }#1}
\theoremstyle{named}

\theoremstyle{definition}
\newtheorem{definition}{Definition}

\theoremstyle{remark}
\newtheorem*{remark}{Remark}

\title{Tiling the 4-ball with knotted surfaces}
\author{James Ross, Hannah Schwartz, Andrew Ye }
\date{}

\begin{document}

\maketitle

\begin{abstract}
We show that for any closed, orientable surface $K$ smoothly embedded in $\mathbb{R}^4$, the unit $4$-ball $B^4 \subset \mathbb{R}^4$ can be tiled using $n \geq 3$ tiles each congruent to a regular neighborhood (with corners) of a surface smoothly isotopic to $K$. This gives a 4-dimensional analog of tilings of the $3$-ball that were constructed in the 90's using congruent knotted tori. 
\end{abstract}

\section{Introduction}

Tilings of the plane have been studied extensively both academically and recreationally; see \cite{gardner} and \cite{thurston} for classic examples. In contrast, topologically interesting tilings of Euclidean $3$-space appeared in the literature starting around the mid-90's, when Adams \cite{Adams, Adams2}, Schmitt \cite{Schmitt, schmitt2} and Kuperberg \cite{Kuperberg} began independently constructing explicit tilings of the $3$-dimensional cube using congruent knotted solid tori. 

Soon after, Oh \cite{Oh} proved the general result that \emph{the $3$-dimensional cube can be tiled using $n$ congruent neighborhoods of any knot type, for any $n \geq 3$.} Such a tiling is called monohedral; see Definition \ref{monohedraltile}. Now we prove the analogous result in $4$-dimensions. 

\begin{theorem}[4D Tiling Theorem] 
Let \(F\) be any closed, orientable surface smoothly embedded in $\mathbb{R}^4$ considered up to smooth isotopy. For any $n\geq 3$, the unit $4$-ball in $\mathbb{R}^4$ admits a monohedral tiling with $n$ tiles, each of which is congruent in $\mathbb{R}^4$ to a neighborhood of the surface $F$.
\label{4dmain}
\end{theorem}

In particular, this answers Question 6 posed by Adams \cite{Adams}, asking for higher-dimensional analogues of his decompositions of cubes into knotted solid tori, and in particular, if the $4$-dimensional hypercube could be decomposed into
thickened knotted $2$-spheres. Up until now, to the authors' knowledge, this had only been partially addressed by Oh \cite[Theorem 3.1]{Oh} in the case where $F$ is a twist-spun knot.

Our proof relies on the existence of bridge trisections due to Meier and Zupan \cite{Other-Meier}. Such a trisection gives a way to decompose any closed surface smoothly embedded in the $4$-ball. For this reason, we must begin with a \emph{smooth} surface, whereas the tiles we recover are \emph{smooth manifolds with corners}. The crux of our argument relies on the following result.

\begin{theorem}[Complement Covering Theorem] 
Let \(F\) be any closed, orientable surface smoothly embedded in $\mathbb{R}^4$ considered up to smooth isotopy. Any smoothly embedded $4$-ball with corners in $\mathbb{R}^4$ can be decomposed into the union of three tiles: two $4$-balls and a neighborhood of the smooth surface $F$.
\label{4dpart2}
\end{theorem}

Note that ``ball coverings" of compact manifolds have been studied since the $70$'s: for instance both \cite{KT76} and \cite[p.~51]{Adams} provide methods to decompose the complement $B^4 - F$ into three $4$-balls. However, utilizing the structure of a bridge trisection of the surface $F$, we are able to build the complement of a neighborhood of $F$ using only \emph{two} $4$-balls, by choosing a special neighborhood that extends all the way out to the boundary of $B^4$ (unlike the previous constructions). Indeed, this allows us to achieve our lower bound of \emph{three} congruent tiles in Theorem \ref{4dmain}. 

Our work was motivated by other recent results on tilings in dimensions $3$ and higher, particularly by \cite{blair2024reptiles}, which gives an isotopy classification of rep-tiles in all dimensions, generalizing the main result of \cite{3drep}. In addition, our proofs of both of the main theorems above were largely inspired by those in Oh \cite{Oh}. We sketch Oh's argument from our perspective in Section \ref{3dtiling}, generalizing it slightly to show the following.

\begin{theorem}[3D Tiling Theorem]
     Let \(T\) be any trivial tangle\footnote{While working in dimension $3$, we place less emphasis on the category. Our proof will work in either the topological category, or as it does in our 4D theorem: taking the tangle $T$ up to smooth isotopy, and the tiles with corners.} properly embedded in the $3$-ball $B$ equal to $D^2 \times I \subset \mathbb{R}^3$, where $D^2 \subset \mathbb{R}^2$ is the closed unit disk. After finitely many stabilizations of $T$, the $3$-ball $B$ can be tiled using any number $n \geq 3$ of congruent copies of a  regular neighborhood of $T$.
     \label{mainthm}
\end{theorem}

The paper is organized as follows: In Section \ref{background} we introduce necessary terminology. In Section \ref{3dtiling}, we discuss tilings of the $3$-ball, proving the 3D Tiling Theorem and detailing Oh's work from \cite{Oh}. We then extend this construction to tile the 4-ball in Section \ref{4d} and prove our two main Theorems \ref{4dmain} and \ref{4dpart2}.

\smallskip 

\noindent
{\bf Acknowledgments.} All four authors would like to thank the North Carolina School of Science and Mathematics, where they conducted this research, as well as our Dean of Mathematics Beth Bumgardner, who champions both her faculty and students to pursue higher level research in public high school, and Manya Nallagangu for her work on the initial stages of this project. The third author HS thanks Ryan Blair, Alexandra Kjuchukova, and Patricia Cahn for introducing her to many interesting open questions about tilings. During this project, she was supported by NSF grant DMS-1502525.

\section{Tilings, Tangles, and Trisections} \label{background}

We work in a variety of categories; as mentioned in the introduction, this is forced upon us by the fact that we utilize smooth techniques, yet produce tiles with corners. Unless otherwise noted, all surfaces and isotopies will be smooth, and all tiles will be smooth at all points except at the ``corners", where they are locally modelled by spaces $[0, \infty)^k \times \mathbb{R}^{n-k}$ with $n$ the dimension of the tile and $0 \leq k \leq n$. We begin by discussing the general notion of a tiling, before specializing to more specific settings.  

\begin{figure}
    \centering
\includegraphics[width=0.6\linewidth]{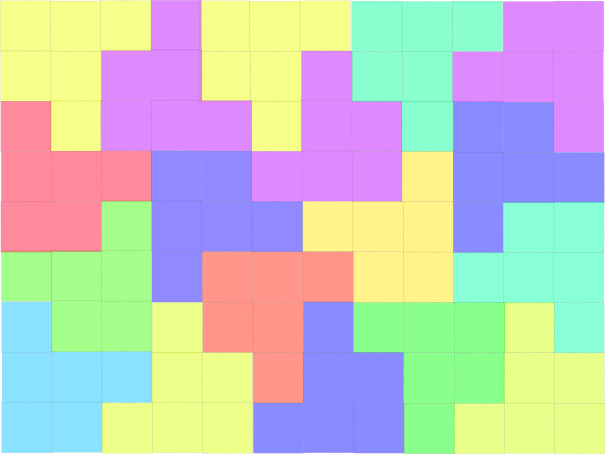}
    \caption{A \textbf{monohedral tiling}, see Definition \ref{monohedraltile}, of a rectangle in $\mathbb{R}^2$, where all tiles are related by an isometry of the plane (rotation, translation, reflection).}
     \label{tilingpic}
\end{figure}

\begin{definition}
A {\bf tiling} of a set $X \subset \mathbb{R}^n$ is a decomposition of $X$ into a  union of at most countably many compact, connected subsets $T_1, T_2, \dots \subset X$ with non-overlapping, non-empty interiors. Our tilings will always be finite. Each such subset $T_i$ is called a {\bf tile}. If all the tiles $T_i$ are congruent to a single {\bf prototile} $T$ via some isometry of $\mathbb{R}^n$, then the tiling is referred to as a {\bf monohedral tiling}. See the 2D tiling in Figure \ref{tilingpic} for an example. 
\label{monohedraltile}
\end{definition}

Our ultimate goal will be to construct tiles congruent (via some isometry of the ambient Euclidean space) to regular neighborhoods of either knotted curves in the $3$-ball, or surfaces in the $4$-ball. Therefore, our final $3$-dimensional tiles will be homeomorphic to a solid torus $S^1 \times D^2$, while our final $4$-dimensional tiles will be homeomorphic to $4$-manifolds of the form $F^2 \times D^2$, 
where $F$ is a closed, orientable surface.

We begin by focusing on how to construct tilings of the $3$-ball, for any given tangle. We will detail this construction in the next section, but develop all necessary terminology here. For convenience below, we fix a homeomorphism $B^3 \simeq D^2 \times [0,1]$ as in the statement of the 3D Tiling Theorem, and a height function on $B^3$ given by the $I$ factor. 

\begin{definition}
    A {\bf tangle} $T$ is a union of $n>0$ arcs properly embedded in the $3$-ball, i.e. so that $\partial T \subset \partial B^3$. A tangle is {\bf trivial} when each arc has a single maximum with respect to the height function on $B^3$. See Adams \cite{Adams3} for a more extensive introduction to tangles and knots.
    \label{tangle}
    \end{definition}

    We will assume not only that all of our tangles are trivial, but also that the maxima of the tangle occur at the same height (this can always be arranged via an isotopy of the tangle rel boundary). We will also consider only projections in which at most one crossing occurs at each height. 
    
    Such a projection of a tangle can be decomposed into the union of three types of ``sections". We divide a diagram of a trivial tangle $T$ into sections as illustrated in Figure \ref{stabletangle}. Analogous sections are described by Oh in \cite{Oh}.  

\begin{definition} \label{sections}
    {\bf Trivial sections} of a tangle diagram consist of $2n$ parallel strands, labeled (c) in our figure, one pair for each arc of the tangle. The {\bf crossing sections} are regions with $2n$ strands where two adjacent strands cross, labeled (b), and the {\bf top section} is where the $n$ maxima of the $n$ arcs occur, labeled (a). \end{definition} 

\begin{figure}[H] 
    \centering 
    \includegraphics[scale= 0.7]{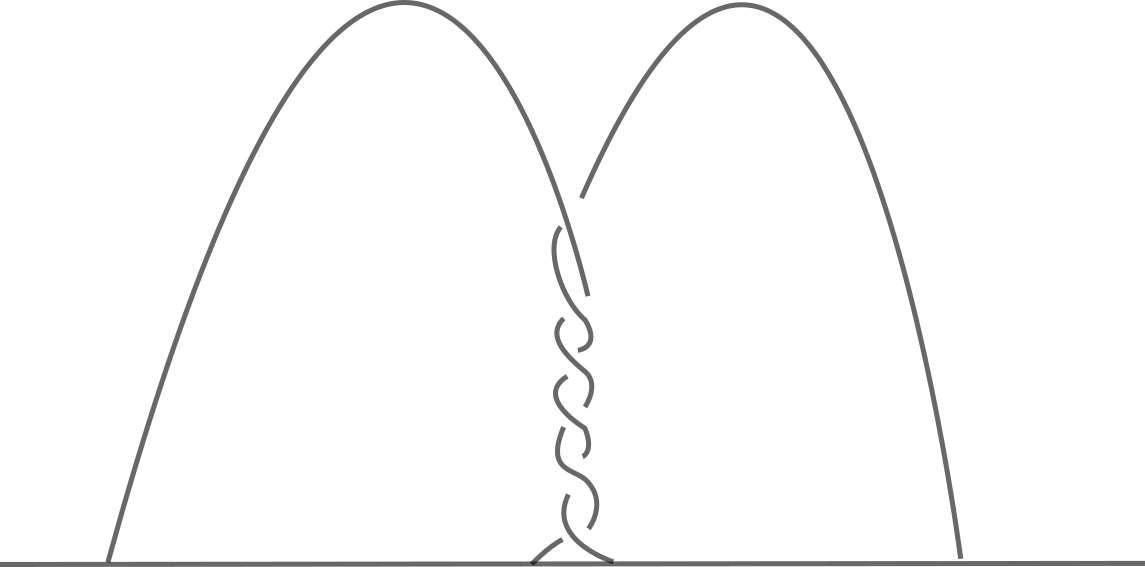}
    \caption{An example of a trivial tangle with two arcs and four strands.}
    \label{tangle1}
    \end{figure}

\begin{definition}
\label{stab}
A diagram of a tangle $T$ is {\bf stabilized} by adding finitely many extra crossing-less single arcs to the diagram, as shown in Figure \ref{stabletangle}. We will refer to the strands added by the stabilization as {\bf helper strands} since they will play a big role in the construction of our tiling.  
\end{definition}

\begin{figure}[H]   
        \centering
        \includegraphics[scale = 0.7]{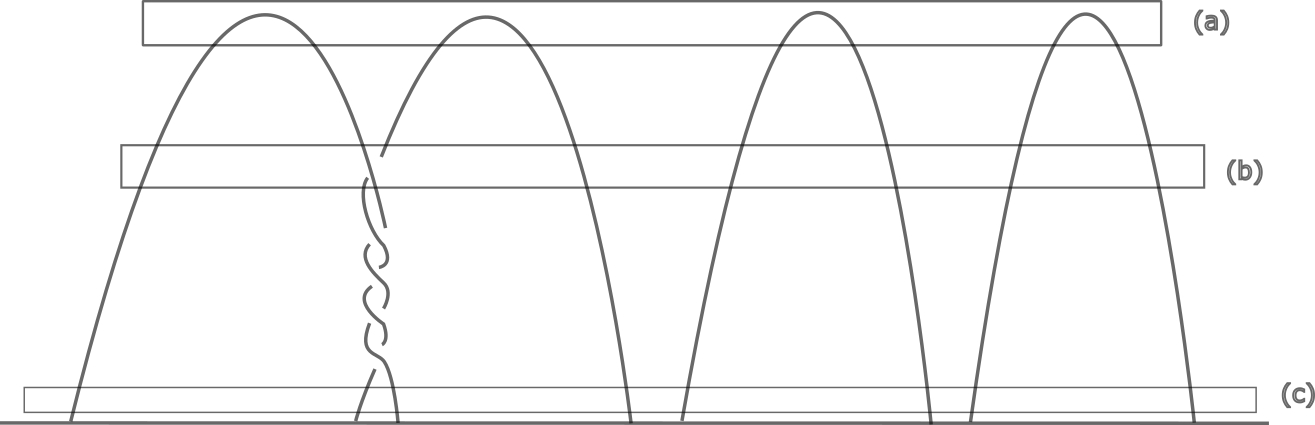}
        \caption{A  tangle stabilized twice, where section (a) is the top section, section (b) is a crossing section, and section (c) is a trivial section.}
        \label{stabletangle}
\end{figure}

It is well known that any knot in $B^3$ can be decomposed into the union of two trivial tangles; bridge splittings of links were first introduced by Schubert \cite{Schubert} in the 50's, and have been extensively studied since. However, only relatively recently has an analog of this statement been given for knotted surfaces. In particular, Meier and Zupan \cite{Other-Meier} showed in 2015 that under certain conditions, a collection of \emph{three} trivial tangles specify a closed surface smoothly embedded in the $4$-ball. 

Such a collection of tangles, along with the decomposition of the surface that it induces, is called a ``bridge trisection" of the surface. This construction is a natural extension of the method of trisecting smooth $4$-manifolds, first defined in the closed case by Gay and Kirby \cite{Gay} in 2016 and later generalized to the compact case by Castro, Gay, and Pinz\'on-Caicedo \cite{reltri}. We utilize the simplest trisection of the $4$-ball, defined in detail below.

\begin{definition}\label{tridef}
    The \textbf{$0$-trisection of $B^4$} consists of a collection of three 4-balls $X_0$, $X_1$ and $X_2$ disjoint in their interiors such that: 
    \begin{enumerate}
        \item The union $X_0 \cup X_1 \cup X_2$ is equal to $B^4$.
        \item Each intersection $  P_{ij} = \partial X_i \cap \partial X_j$ is homeomorphic to a $3$-ball.
        \item  The triple intersection $\Sigma = P_{01} \cap P_{12} \cap P_{20}$ is a disk. 
    \end{enumerate}
We refer to each $3$-ball $P_{ij}$ as a {\bf page} of the trisection, and to the disk $\Sigma$ of triple intersection as the {\bf binding}. This is in reference to the open book decomposition of $B^4$ induced by the trisection, illustrated in Figure \ref{mainthm}.  
\end{definition}

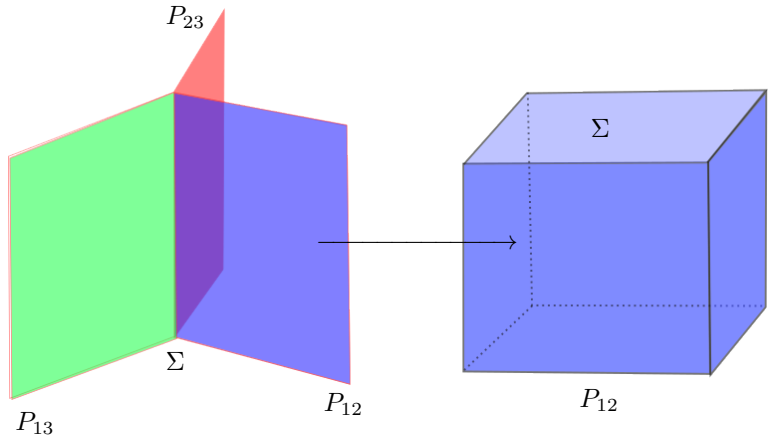
\begin{figure}
\label{zerotripic}
\begin{center} \includegraphics[scale=0.5]{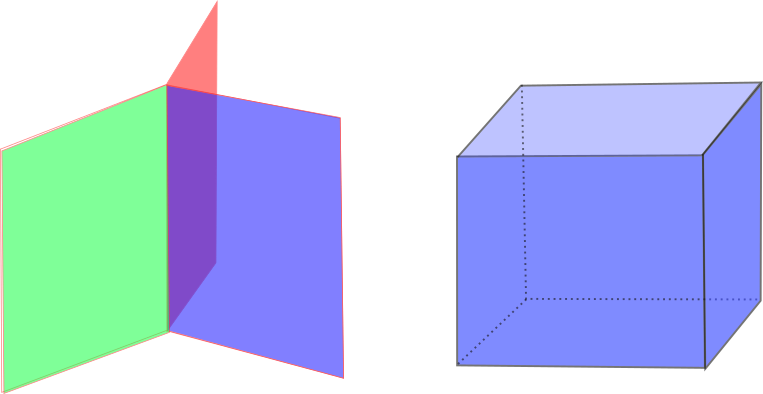}
\put(-6,2){\(\xrightarrow{\hspace{7em}}\)}
\put(-8,0.4){\(\Sigma\)}
\put(-5.9,-0.15){\(P_{12}\)}
\put(-10,-0.4){\(P_{13}\)}
\put(-8,5){\(P_{23}\)}
\put(-2.5,-0.1){\(P_{12}\)}
\put(-2.35,3.5){\(\Sigma\)}
\end{center}
\caption{A PL depiction (drawn down one dimension) of the $0$-trisection of a \(4\)-ball into three separate \(4\)-balls $P_{ij}$, whose triple intersection is the binding \(\Sigma\)} of the open book decomposition shown on the left.
\end{figure}

To accompany this usual topological definition of the trisection, we addend a geometric version. This will later come in handy when constructing our tiles as explicit subsets of $\mathbb{R}^4$.

\begin{definition} \label{defgeostd}
    The \textbf{geometrically standard $0$-trisection of the unit ball $B^4 \subset \mathbb{R}^4$} is a decomposition of $B^4$ into three 4-balls $X_0$, $X_1$ and $X_2$, where
         $$X_i = \{(r, \theta, z,t) \in B^4  ~|~  \frac{2\pi}{3}i \leq \theta \leq  \frac{2\pi}{3}(i+1)\}$$ 
         such that $i+1$ is taken mod $3$. Notice that $X_0, X_1$ and $X_2$ satisfy the conditions of Definition \ref{tridef}, and therefore give a trisection of the $4$-ball. 
\end{definition}

We end this section by introducing the notion of a bridge trisection, the 4-dimensional analog to a bridge splitting of a link, first introduced by Meier and Zupan in \cite{Other-Meier}. 

\begin{definition} \label{btri}
    A \textbf{bridge trisection} $(T_{01}, T_{12}, T_{20})$ of a closed surface $F$ smoothly embedded in the $4$-ball consists of a collection of three properly embedded tangles $T_{ij} \subset  P_{ij}$ such that
    \begin{enumerate}
        \item Each tangle $T_{ij}$ is trivial, as in Definition $2.2$.
        \item The union $L_j = T_{ij} \cup -T_{jk}$ is an unlink in the $3$-ball $P_{ij} \cup_\Sigma -P_{jk}$, which bounds a collection \(\mathcal{D}_j\) of disjointly embedded spanning disks.
        \item The surface $F$ is isotopic to the union $\mathcal{D}^+_0 \cup \mathcal{D}^+_1 \cup \mathcal{D}^+_2$, where each $\mathcal{D}^+_j$ is a collection of disks properly embedded the $4$-ball $X_j$ simultaneously isotopic rel boundary to the spanning disks $\mathcal{D}_j $. 
    \end{enumerate}
\end{definition}

Heuristically, bridge trisection diagrams aid us in visualizing and deconstructing embedded surfaces in the $4$-ball into more manageable pieces. In \cite{Other-Meier}, Meier and Zupan show that any closed surface smoothly embedded in the $4$-ball admits a trisection diagram (albeit non-uniquely). We will use this fact, and the notation above, throughout the remaining sections.

\section{Tiling the 3-ball with trivial tangles} \label{3dtiling}

We begin by presenting our proof of Theorem \ref{mainthm}, our slight generalization of Oh's result from \cite{Oh}.  Many ideas from our interpretation of his argument are utilized in our argument one dimension higher.

\smallskip
\noindent
\emph{Proof of the 3D Tiling Theorem.} Start by fixing a diagram of $T$, and stabilize this diagram $n$ times, so that the tangle $T$ now has $2n$ arcs. We consider all endpoints of the tangle in $D^2 \times \{0\} \subset \partial B$. We shall use a construction based on Oh's from  \cite{Oh} to show that any number of congruent copies of a neighborhood of $T$ (the stabilized version) can tile this \(3\)-ball in $\mathbb{R}^3$. Our tiling will be composed of smaller tilings of cylindrical $3$-balls glued together appropriately; we will refer to each of these smaller tilings as \emph{sections} of the main one. Each such section will correspond to a section of our fixed projection of the tangle $T$, via the construction detailed below. Although we begin by building only $3$ congruent tiles, we will later explain how an analogous tiling can be made using any number $n \geq 3$ of tiles.  \\

\noindent \emph{Constructing a \textbf{trivial section} of the tiling:}  

\vspace{.3cm}

\begin{figure}[H] 
\centering
\begin{minipage}{.45\textwidth}
    \centering
\includegraphics[width=0.55\textwidth]{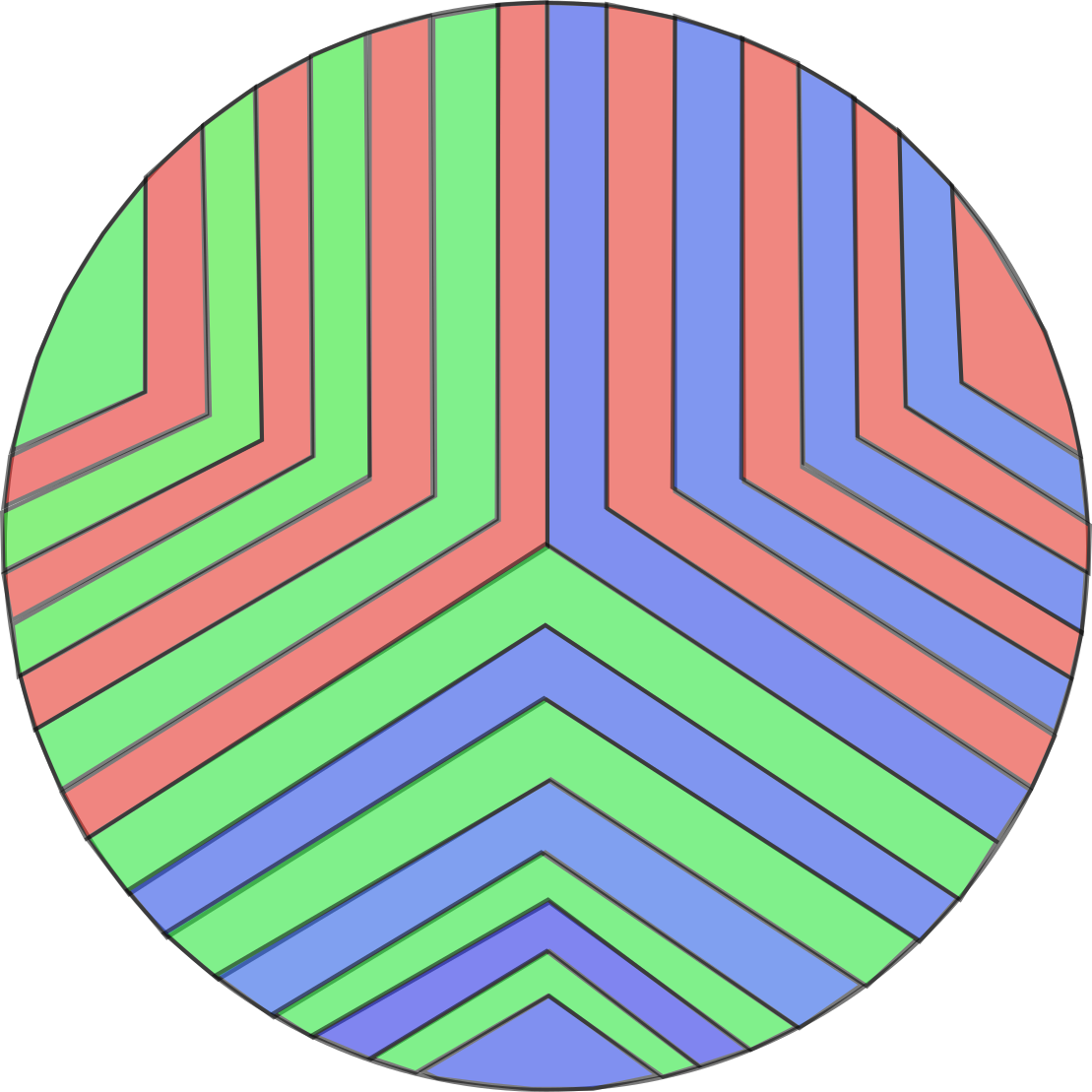}
    \caption{Top down view of a \emph{trivial section}.}
    \label{toptile}
\end{minipage}
~~~
\begin{minipage}{0.45\textwidth}
    \centering
\includegraphics[width=0.795\textwidth]{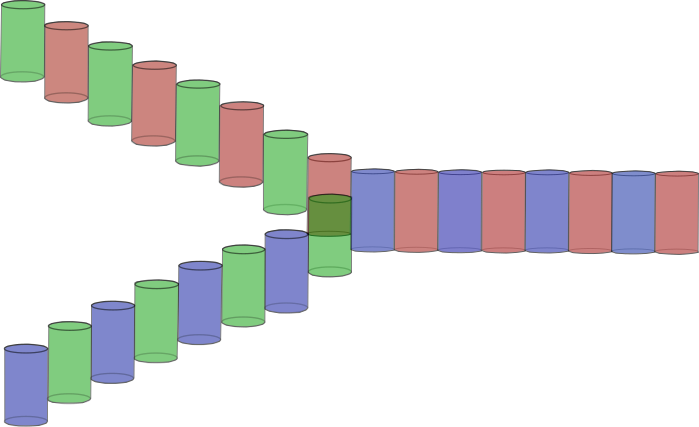}
    \caption{Side view of a \emph{trivial section}.}
    \label{fig2}
\end{minipage}
\end{figure}

Figures  \ref{fig2} and \ref{toptile} illustrate what we call a \emph{trivial section} of our final tiling. In particular, Figure  \ref{fig2} gives a ``side view" of a trivial section,  whereas Figure \ref{toptile} depicts the same tiling (up to homeomorphism) from the ``top-down".

We use color to partition our decomposition of the $3$-ball into three congruent tiles, each of which is the union of a collection of cylindrical columns. In turn, each column corresponds to one strand in a trivial section of the tangle $T$. For instance, in the examples below, the eight green cylinders correspond to the eight strands of a single copy of a trivial section of $T$. Four of these cylinders (those in one ``branch" of the {\bf Y} from Figure \ref{fig2} or one ``slice" of the disk in Figure \ref{toptile}) correspond to the original strands of $T$, while the other four are the helper strands, as described in Definition \ref{stab}. 

\noindent \emph{Constructing a \textbf{crossing section} of the tiling:}

\vspace{.3cm}

\begin{figure}[H] 
\centering
\begin{minipage}{.45\textwidth}
    \centering
    \includegraphics[width=0.55\textwidth]{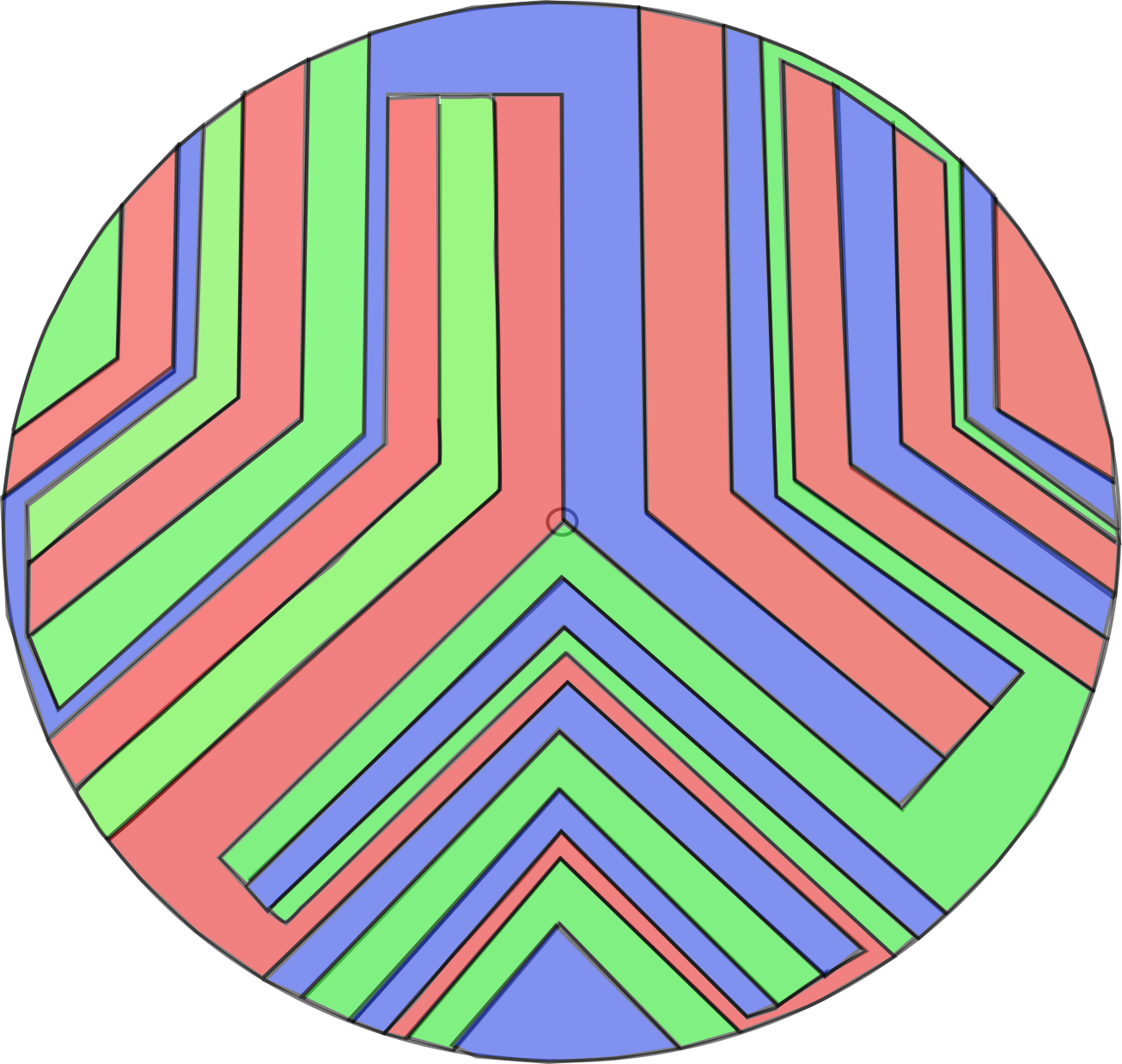}
    \caption{A top down view of the first and third levels of the \emph{crossing section}.}
    \label{topdowncross}
\end{minipage}
~~~
\begin{minipage}{0.45\textwidth}
    \centering
    \includegraphics[scale = 0.18]{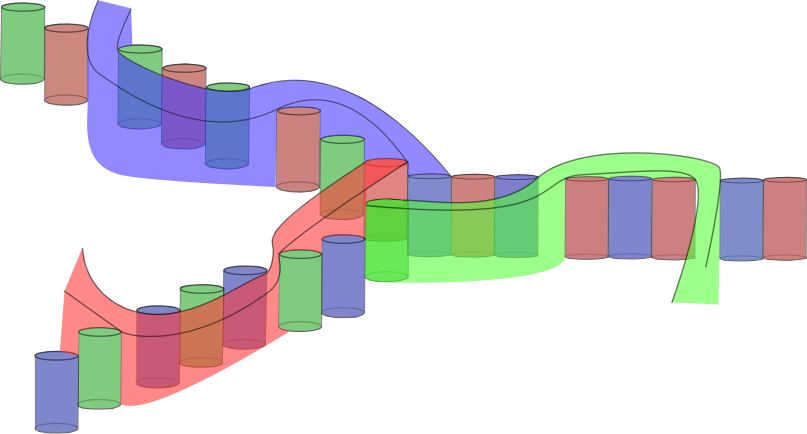}
    \caption{A side view of the first and third levels of the \emph{crossing section}.}
    \label{sideclose}
\end{minipage}
\end{figure}

\begin{figure}[H]
    \centering
    \includegraphics[scale = 0.5]{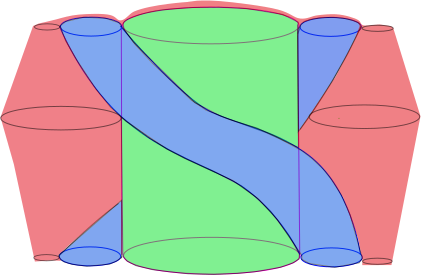}
    \caption{A side view of one crossing in the second level of the tiling of the \emph{crossing section of $T$}.}
    \label{closecrossing}
\end{figure}

We construct this new section by stacking three \emph{levels} of smaller tilings of $B^3$. The first and third levels are identical: this tiling of $B^3$ by parallel cylinders is shown in Figures \ref{topdowncross} and \ref{sideclose}. Essentially, the first level arranges the tiles into position to prepare for the second level where the ``crossing" actually occurs, as pictured in Figure \ref{closecrossing}. Tiles of each color are required to prevent tiles of the same color from touching while realizing a crossing. Note that this is what necessitates our use of helper strands, and at least three distinct tiles. The third level of the crossing section is identical to the first. Stacking these levels together sequentially produces our tiling for the crossing section, in which each cylindrical component corresponds to one strand of a crossing section of the projection of $T$.


In Figure \ref{closecrossing}, the tiling in the crossing section corresponding to a single crossing of $T$ is shown. The helper strands (copies of whose neighborhoods are shown in red and green) isolate the neighborhoods of the blue strands involved in the crossing. Up to permuting colors, this illustrates each crossing that occurs in the second \emph{level} of the crossing section.

\noindent \emph{Constructing the \textbf{top section} of the tiling:} 

\vspace{.3cm}

\begin{figure}[H] 
\centering
\begin{minipage}{.45\textwidth}
    \centering
    \includegraphics[width=0.55\textwidth]{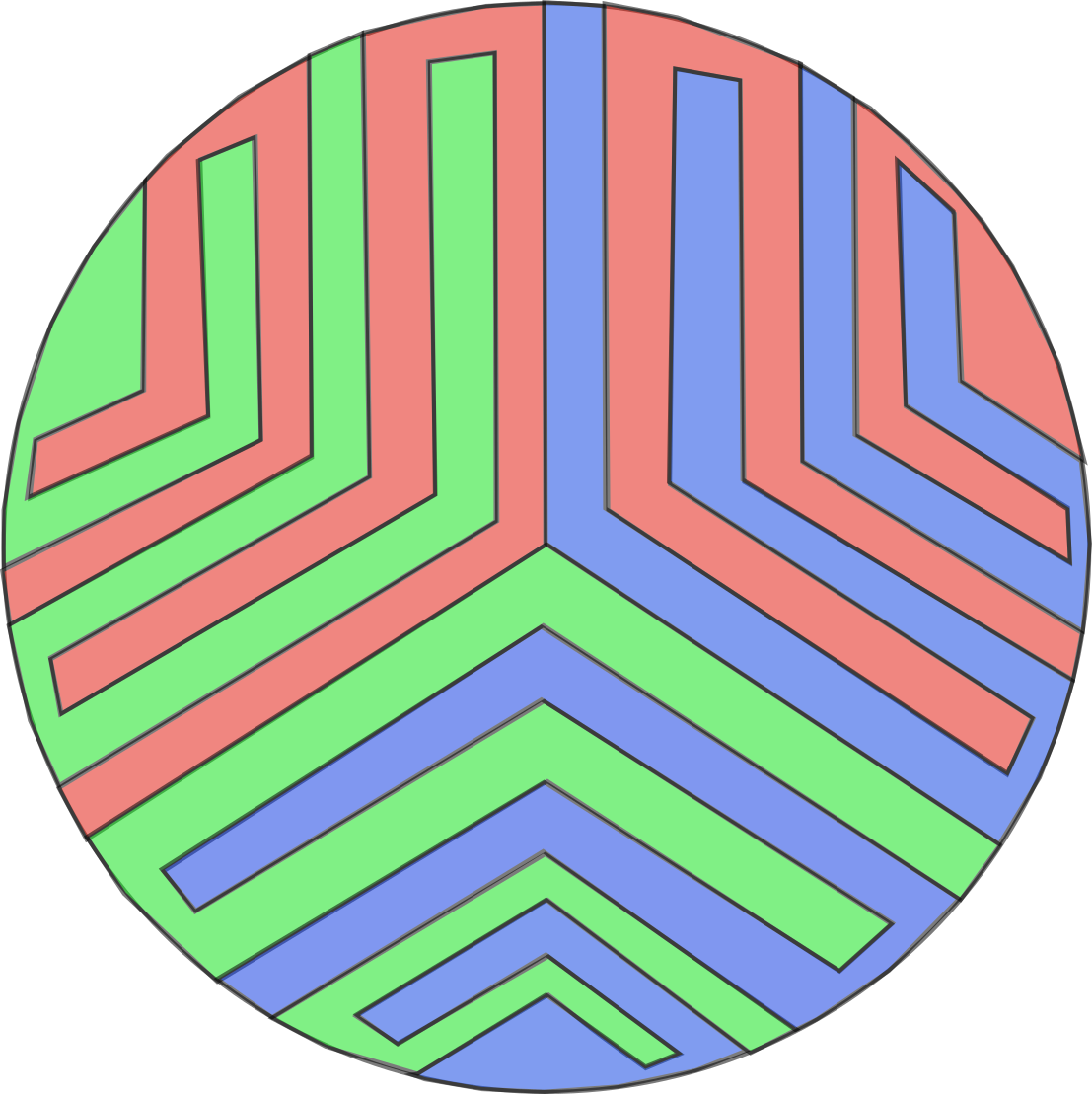}
    \caption{A top down view of the tiling that forms the \emph{top section}.}
    \label{topdownclose}
\end{minipage}
~~~
\begin{minipage}{0.45\textwidth}
    \centering
    \includegraphics[scale = 0.235]{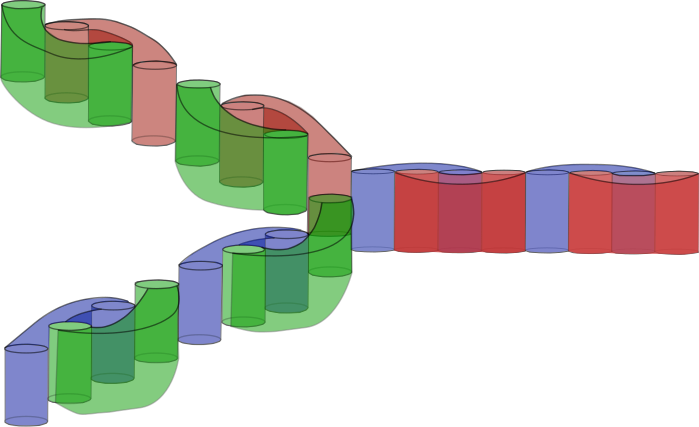}
    \caption{A side view of the tiling that forms the \emph{top section}.}
    \label{sidecross}
\end{minipage}
\end{figure}

The \emph{top section} of our tiling of $B^3$ is constructed in order to connect pairs of  components of the same color which are disconnected in the first two sections. Each tile in this section will be identified with a neighborhood of the maxima of the tangle $T$. A top-down and side view of the top section is illustrated in Figures \ref{topdownclose} and \ref{sidecross}. 

\noindent \emph{\textbf{``Stacking'' the Tilings:} Stack the tilings of $D^2 \times I$ from Steps $1$-$3$.} 
\vspace{.3cm}

A finite sequence of the tilings of $D^2 \times I$ created in the previous steps can be stacked so that the sections used as tiles mirrors the order of these sections in the tangle. By construction, consecutive levels can be placed so that tiles of like colors overlap on the boundary. This produces our final tiling of $D^2 \times I \simeq B^3$ using $3$ congruent copies of a neighborhood of the tangle $T$. 

\vspace{3mm}

This construction can be generalized to tile the disk using $n\geq 3$ congruent tiles, by adding extra ``branches" or ``slices" to each section of the tiling. Refer to Figure 5 of Oh \cite{Oh} for more detail.

\begin{remark}
    The tiling produced in the previous theorem depends not only on the tangle $T$, but also on our initial choice of projection. 
\end{remark}

As mentioned, the argument given above largely follows the construction given in Oh \cite{Oh}. As such, it is not surprising that our generalization implies Oh's result as a corollary. 

\begin{corollary} \label{tilingwithL}
        There exists a tiling of the $3$-ball by \(n\) congruent copies of any link  \(L \subset B^3\). 
\end{corollary}

\proof 
Take a bridge splitting of the tangle $L$ into two trivial tangles $T_1$ and $T_2$; that this can be done is a classical result due to \cite{Schubert}. Stabilizing each of these tangles sufficiently, the $3$-ball on either side of the bridge sphere can be tiled using $n$ tiles, all of which are congruent to a neighborhood of the (stabilized) tangles $T_1$ and $T_2$, by Theorem \ref{mainthm}. Identifying corresponding pairs of tiles from each $3$-ball along their boundaries gives a tiling of the $3$ ball, in which each of the $n$ resulting tiles is congruent to some neighborhood of $L$. 

\qed

\section{Tiling the 4-ball with knotted surfaces}
\label {4d}

We now prove our 4-dimensional analog of the results from the previous section. To prepare for our proofs of both theorems, let $F$ be any closed, orientable surface smoothly embedded in $\mathbb{R}^4$, considered up to smooth isotopy. Equip the unit $4$-ball in $\mathbb{R}^4$ with the geometrically standard $0$-trisection as in Definition \ref{defgeostd} in which all $X_i$ are congruent. Then smoothly isotop $F$ into the unit $4$-ball so that it lies in bridge position with respect to this trisection, with bridge trisection $(T_{01}, T_{12}, T_{20})$ as in Definition \ref{btri}. We will work with a fixed diagram of each tangle $T_{ij}$.\\

The final tiling we construct to prove the 4D Tiling Theorem (see Figures \ref{tilingtri} and \ref{pinwheelpic}) will be a ``pinwheel" akin to those in \cite{gpinwheel} and \cite{phpinwheel} used to decompose the $4$-ball into multiple identical pieces. In particular, we follow the general steps outlined below: 

\begin{enumerate}
    \item Decompose the $4$-ball into the union of one interesting piece $P$ and some number $n-1$ of smaller $4$-balls $B_1, \dots, B_{n-1}$.
    \item Take (at least) $n$ copies of the $4$-ball, each with such a decomposition, thinking of each copy as the ``spokes" of the pinwheel around a central hub.  
    \item The boundary connected sum of the spokes can be decomposed into $n$ identical tiles, each of which is the boundary connected sum of the piece $P$ in one spoke, together with copies of $B_1, \dots, B_n$ from the remaining spokes. 
\end{enumerate}

A similar construction, specifically of a tiling, was also suggested by Adams in \cite[p.~51]{Adams}. Central to this strategy is the decomposition of the ball in Step 1, which we construct below.
\\

\noindent \emph{Proof of the Complement Covering Theorem}: Note that it is sufficient to produce such a decomposition of any \emph{one} embedded $4$-ball with corners in $\mathbb{R}^4$. For, Palais' disk theorem \cite{palais} can be used to produce an isotopy taking this $4$-ball to any other\footnote{Pre- and post-composing with an isotopy that smooths and then ``un"-smooths the corners, each supported on a collar of the boundary away from the surface $F$.} restricting to a smooth isotopy on the surface $F$. This isotopy carries the decomposition of the original $4$-ball to a corresponding collection of tiles for the second $4$-ball.

\begin{figure}[ht] 
\centering
\includegraphics[scale=0.5]{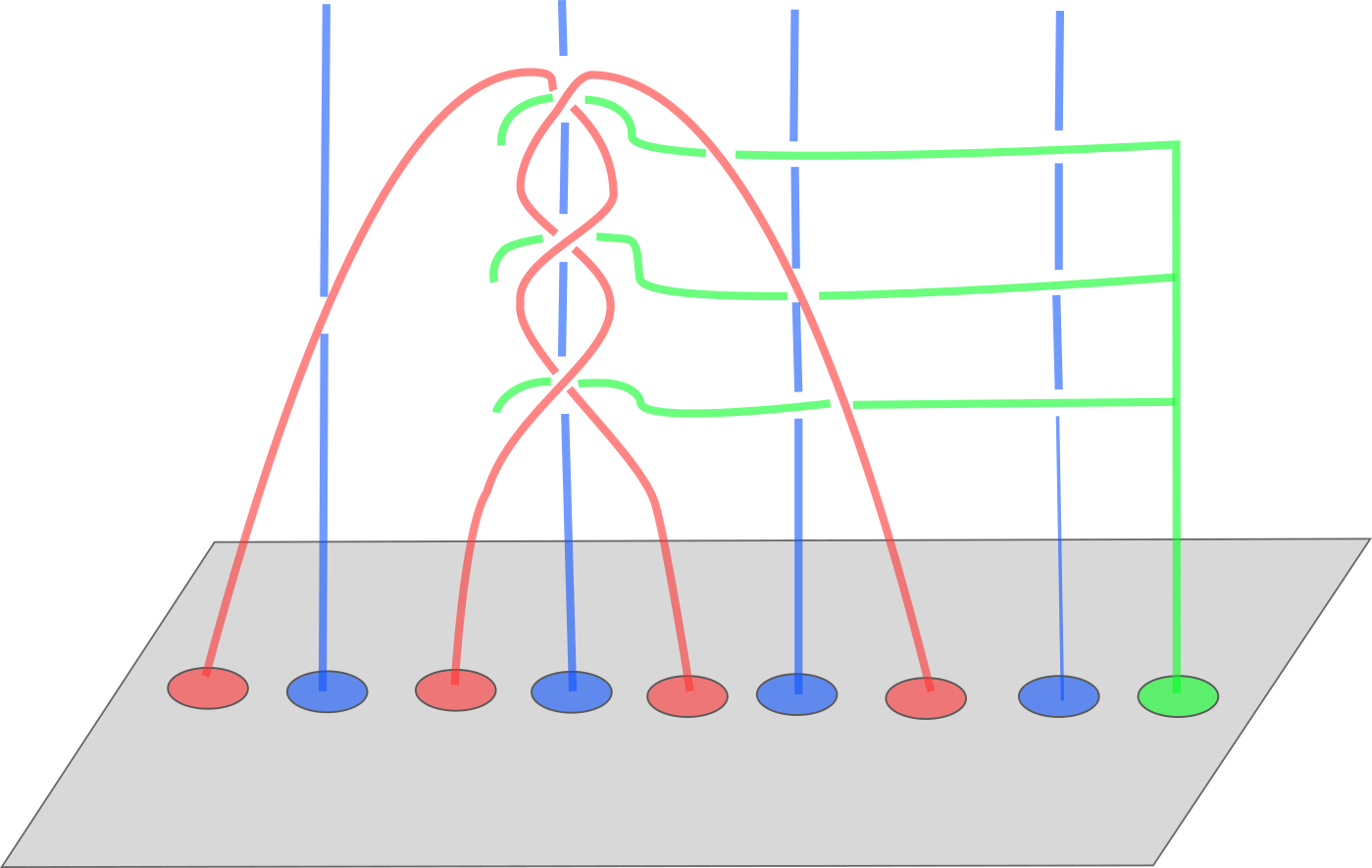}
\put(-8.4,2.5){\textcolor{red}{Tangle}}
\put(-7.9,5){\textcolor{blue}{Loom}}
\put(-1, 4){\textcolor{green}{Scaffolding}}
\put(-1, 0){$\Sigma$}
\caption{A depiction of the tangle $T_{ij}$ (in red), the loom (in blue) and the scaffolding (in green) in the page $P_{ij}$. Our final tiling of $P_{ij}$ is constructed by taking sufficiently large neighborhoods of each of these pieces.}
\label{pageyloom}
\end{figure}

\begin{figure}[ht] 
\centering
\includegraphics[scale=1]{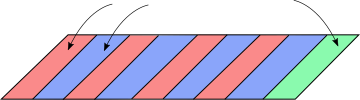}
\put(-6.4,2.6){\textcolor{red}{Red}}
\put(-5.45,2.5){\textcolor{blue}{Blue}}
\put(-2.8, 2.6){\textcolor{green}{Green}}
\put(-1, 0){$\Sigma$}
\caption{Our final tiling of each page $P_{ij}$ restricted to binding $\Sigma$. Note that this tiling of $\Sigma$ is independent of the page $P_{ij}$.}
\label{boundarytile}
\end{figure}
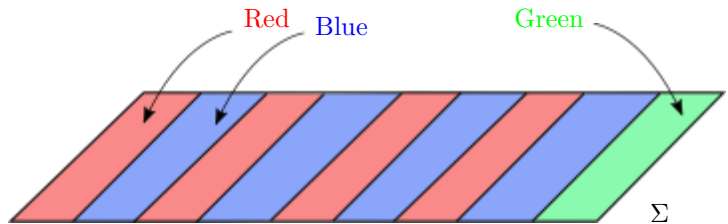

The tiling we produce will be strictly contained in the unit $4$-ball. At each step in the construction, we will use three colors (red, green, and blue) to keep track of which regions will contribute to which of the three tiles in our final decomposition. We begin by tiling each page \(P_{ij}\) of the geometrically standard trisection of $B^4$. 

First, color the tangle $T_{ij}$ in one color (say, red). Next, construct a ``loom'' in blue, consisting of the union of parallel, vertical line segments aligned on a plane behind the red tangle diagram. Position these blue line segments so that each point of their intersection with the binding lies between two consecutive endpoints of the red tangle.  Finally, construct a green ``scaffolding'' which will serve a purpose analogous to the ``helper strands'' from the proof of \thref{mainthm}. The scaffolding is the union of a single vertical strand intersecting the binding once, and a collection of horizontal strands each of which wraps around a single crossing of the tangle \(T_{ij}\). See Figure \ref{pageyloom} for an explicit example.

 Taking sufficiently large closed neighborhoods (with corners) around the tangle, the loom, and the scaffolding gives a tiling of the $3$-ball $P_{ij}$ by a collection of multiple disjoint red, blue, and green $3$-balls. Observe that the red tile is a neighborhood of the tangle $T_{ij}$; likewise the blue tile is equal to a neighborhood of the loom, and the green tile a neighborhood of the scaffolding. In order for us to piece these tilings together, we arrange that the tilings of the binding \(\Sigma\) induced by the tilings of each page $P_{ij}$ are identical, and do not depend on the tangle specifically. We can achieve this by choosing our tiles so that their intersection with $\Sigma$ is the tiling consisting of rectangular strips of alternating colors, as shown in Figure \ref{boundarytile}. 

Now each $3$-ball \(\partial X_j\) can be tiled by identifying tiles from \(P_{ij}\) and \(P_{jk}\) along the binding $\Sigma$ so that tiles of the same color ``match up".  After this identification, the red tile is a neighborhood of the unlink $L_j = F \pitchfork \partial X_j$ (see Definition \ref{btri}) while the green and blue tiles are each still homeomorphic to a collection of $3$-balls. 
 
 We proceed to extend this tiling to a tiling of a $4$-ball $X_j'$ properly contained in $X_j$. First, identify $X_j$ with the product $\partial X_j \times [0,1]$ so that $B^3 \times \{0\} = P_{ij} \cup -P_{jk}$. We tile the ``upper half" of this product $\partial X_j \times [0,1/2]$ with $4$-dimensional tiles equal to thickened versions of the $3$-dimensional tiles used to decompose $\partial X_j$. In other words, each 4D tile is simply the product of a 3D tile with the interval $[0,1/2]$. 

  \begin{figure}[H]
     \centering     \includegraphics[scale=0.3]{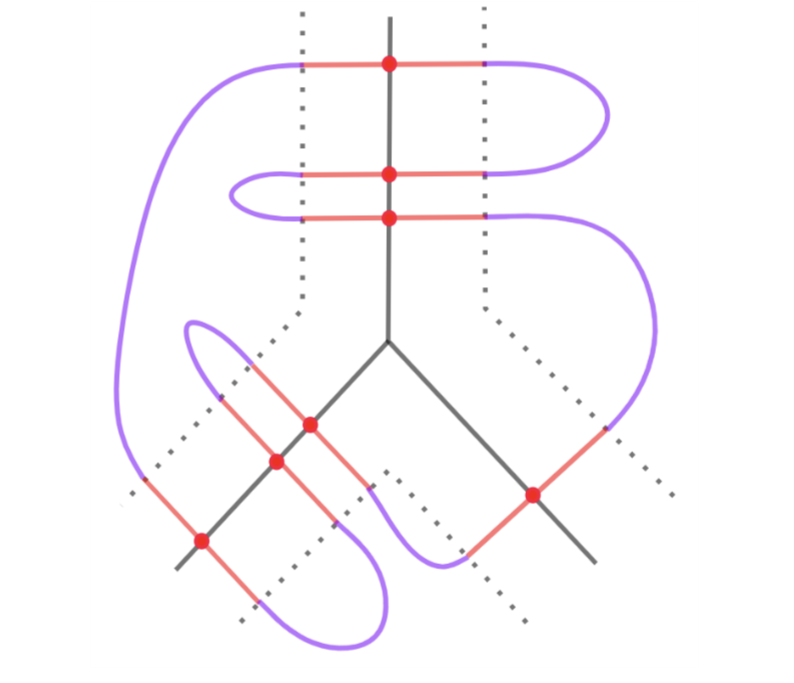}
     \put(-3.3,5.6){\(0\)}
     \put(-4.05,5.65){\(\frac{1}{2}\)}
     \caption{A schematic of the bridge trisection of $F$ used to construct the tiling in our proof of the Complement Covering Theorem. The tangles $(T_{01}, T_{12}, T_{20})$ are shown in red as points on the ``spine" of the trisection at the $0$ level. The disks $\mathcal{D}^+_j$ are illustrated in purple, below the $\frac{1}{2}$ level. }
     \label{fig:enter-label}
 \end{figure}

Now, recall from Definition \ref{btri} that  $F \cap X_j $ is a union of disjoint, properly embedded disks $\mathcal{D}^+_j$, that are ``pushed in" copies of embedded disks $\mathcal{D}_j$ bounded by the unlink \(L_j\) in each $3$-ball $\partial X_j$. Let $N(\mathcal{D}_j)$ be a regular neighborhood of these disks in $\partial X_j$, and add the region $N(\mathcal{D}_j) \times [1/2,1]$ to the red tile. After this addition, the red tile in $X_j$ is equal to a neighborhood of a set of disks isotopic to $\mathcal{D}_j^+$ rel boundary. 

By construction, the tilings of each pair $X_i'$ and $X_j'$ are identical along the page $P_{ij}$. Therefore, we can decompose the $4$-ball $X_0' \cup X_1' \cup X_2'$ into one red, one blue, and one green tile, each of which is the identification of the similarly colored pieces from our decompositions of $X_0', X_1'$ and $X_2'$ along their boundaries. The red tile is a neighborhood of the smooth surface $F$, and both the blue tile and the green tile are $4$-balls. 

\qed

\smallskip 

Finally, we proceed to prove our main result. 

\smallskip 

\noindent \emph{Proof of the 4D Tiling Theorem.} Let $X_0 \cup X_1 \cup X_2$ be the geometrically standard $0$-trisection of the unit $4$-ball. Since $X_0$ is a $4$-ball with corners, by the Complement Covering Theorem there is a tiling of $X_0$ such that one tile (the ``red") is a neighborhood of the smooth surface $F$ up to smooth isotopy, and the two remaining tiles (colored ``green" and ``blue") are $4$-balls. Let $T_R$ denote this tiling, to record that the neighborhood of the surface $F$ is the \emph{red} tile. 

Note that by construction, all three tiles in $T_R$ have non-empty intersection with $\partial X_0$. Therefore, via a smooth ambient isotopy of $\mathbb{R}^4$ supported only on a collar of $\partial X_0$ and fixing $X_0$ set-wise, we may arrange for the restriction of the tiling on the boundary of $X_0$ to resemble that shown in Figure \ref{tilingtri}. This standardizes the tiling on the portion of $\partial X_0$ intersecting $\partial X_1 \cup 
 \partial X_2$, and pushes the complexity elsewhere. It will also ensure that our final tiles have the correct homeomorphism type.

The $4$-balls $X_1$ and $X_2$ admit analogous tilings, gotten by rigidly rotating the tiling $T_R$ of $X_0$ around the planar axis of $\mathbb{R}^4$ containing the binding $\Sigma$. Permute the colors in each of these tilings, so that the  neighborhood of $F$ is green in $X_1$ and blue in $X_2$.  We refer to these two new tilings as $T_G$ and $T_B$, respectively.

\begin{figure}[H] 
\centering
\begin{minipage}{.45\textwidth}
    \centering
\includegraphics[width=.95\textwidth]{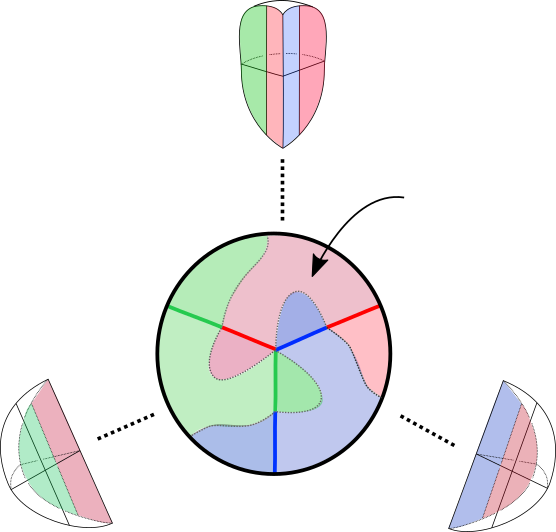}
\put(-3.5,4.25){{$T_R$}}
\put(-.1,0){{$T_B$}}
\put(-1.6,3.4){{Neighborhood}}
\put(-1.1,3){{of $F$}}
\put(-5.5,0){{$T_G$}} 
    \caption{Our final tiling of the unit $4$-ball, the union of the tilings $T_R, T_G$ and $T_B$ of each ``wedge" $X_{i}$. These tilings are chosen to coincide along their boundaries so that each new tile is isotopic to a neighborhood of $F$.}
    \label{tilingtri}
\end{minipage}
~~~
\begin{minipage}{0.45\textwidth}
    \centering
    \includegraphics[scale = 0.21]{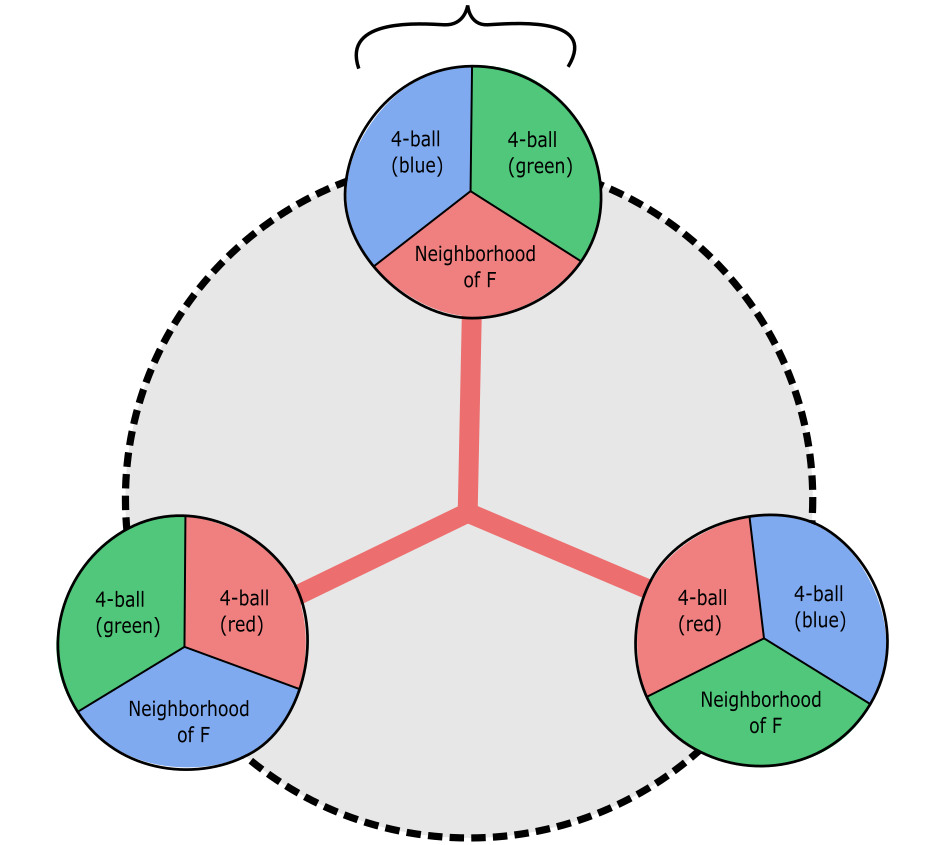}
    \put(-2.8,4.8){{$B^4$}}
    \caption{A schematic of our ``pinwheel" construction of a tiling as a union of identical tilings of $B^4$ boundary summed around a central spoke.}
    \label{pinwheelpic}
\end{minipage}
\end{figure}

The union of the three tilings $T_R$, $T_G$ and $T_B$ (identifying tiles of like colors along their boundaries) gives the final desired tiling of $B^4$ illustrated in Figure \ref{tilingtri}. The result is three congruent tiles, the $i^{th}$ of which is a boundary connected sum of the neighborhood of $F$ in $X_i$ with a $4$-ball in each $X_j$ with $i \not = j$. Note that by using extra colors and inserting additional ``wedges" into the decomposition of the unit $4$-ball, we may extend this construction to build a tiling using any number $n \geq 3$ tiles\footnote{A \emph{minimum} of three distinct tiles is needed, however, in order for the final connected sums to have the correct homeomorphism type.}.\qed

\newpage
\bibliographystyle{plain}
\bibliography{references.bib}

\end{document}